\documentclass{amsart}

\usepackage{amsmath,amssymb,amsfonts,amsthm,graphicx}

\newtheorem{theorem}{Theorem}[section]

\newtheorem{corollary}[theorem]{Corollary}
\numberwithin{equation}{section}

\theoremstyle{remark}

\newcommand{\Ric}{\mathop{\mathrm{Ric}}}

\title{Bounds on the Heat Kernel under the Ricci Flow}

\author{Mihai Bailesteanu}

\begin{document}

\maketitle

\begin{center}
Department of Mathematics, Cornell University, \\ 120 Malott Hall,
Ithaca, NY 14853-4201, USA \\ \texttt{mbailesteanu@math.cornell.edu} \\

\end{center}

\begin{abstract}
We establish an estimate for the fundamental solution of the heat equation on a closed Riemannian manifold $M$ of dimension at least $3$, evolving under the Ricci flow. The estimate depends on some constants arising from a Sobolev imbedding theorem. Considering the case when the scalar curvature is positive throughout the manifold, at any time, we will obtain, as a corollary, a bound similar to the one known for the fixed metric case.  
\end{abstract}

\section{Introduction}

We consider a closed manifold $M$, whose metric evolves by the Ricci flow: \begin{equation}\label{Riccifloweqn}
\frac\partial{\partial t} g(x,t)=-2\Ric(x,t),\qquad x\in M,~t\in[0,T]
\end{equation} 
and we obtain un upper bound for the heat kernel $G(x,t;y,s)$, which is the fundamental solution of the heat equation

\begin{align}\label{heateqn}
\left(\Delta-\frac\partial{\partial t}\right)u(x,t)=0,\qquad x\in M,~t\in[0,T].
\end{align}

Determining bounds for the heat operator on manifolds has been a topic of interest, as it had proven to have many applications. D. Aronson made use of a parabolic Harnack inequality to bound the fundamental solution for a general second-order parabolic operator \cite{DA67}. Later, in their celebrated paper \cite{PLSTY86}, P. Li and S.-T. Yau derived gradient estimates for positive solutions to the heat equation on closed manifolds with bounded Ricci curvature, from which they obtained Harnack inequalities. Further these inequalities were used to get upper and lower bounds on the heat kernel. They considered manifolds with boundaries, satisfying Dirichlet and Neumann boundary conditions, the heat kernel bounds extending to the boundary when the latter was convex. Further, J. Wang derived in \cite{JW97} a global version of gradient estimates when the boundary is nonconvex, and he obtained both upper and lower bounds for the heat kernel satisfying Neumann conditions. 

In geometric analysis, heat kernel estimates, together with Sobolev imbedding theorems, have been proven useful in the study of Ricci flows, especially in the case with surgeries. Since Sobolev imbeddings and inequalities relate the integrability (in some $L^p$ sense) of the derivative of a function to the integrability of the function itself,  they become useful when looking at partial differential equations. They also have proven useful in characterizing the space where the function is defined (for a detailed discussion see, for example, \cite{LSC01}).   

In \cite{CG02}, C. Guenther studied the fundamental solution of the linear parabolic operator $L(u)=(\Delta-\frac{\partial}{\partial t}-h)u$, on compact n-dimensional manifolds with time dependent metric, where $h$ is a smooth space-time function. She proved the uniqueness, positivity, the adjoint property and the semigroup property of this operator, which thus behaves like the usual heat kernel. As a particular case ($h=0$), she obtained the existence and properties of the heat kernel under the Ricci flow.

G. Perelman gave a proof in \cite{GP02} of the pseudolocality theorem, which states that Euclidean looking regions in closed manifolds evolving by Ricci flow remain localized, under some curvature assumptions. In order to prove this, he obtained a differential Li-Yau-Hamilton type inequality for the fundamental solution of the conjugate heat equation $\Delta u+ \frac{\partial}{\partial t}-Ru=0$, where $R$ is the scalar curvature on the manifold. Later S.-Y. Hsu obtained in \cite{SH08}, by a variation of the method introduced by P. Li, S.-T. Yau and J. Wang, a gradient estimate for the solution of the conjugate heat equation on closed manifolds under Ricci flow, and as a consequence, bounds for its fundamental solution. 

Q. Zhang also considered in \cite {QZ06} the conjugate heat equation introduced by Perelman, but after a time reversal: $\Delta u-\frac{\partial}{\partial t}-Ru=0$, and the metric evolving under forward Ricci flow $\frac\partial{\partial t} g(x,t)=2\Ric(x,t)$. He considered a complete manifold, with $\Ric(g(t))\ge k$ and the injectivity radius $i>0$. He obtained un upper bound on the fundamental solution of this equation using the Nash method, without any gradient estimate, and his result depends on the best constants in the Sobolev imbedding theorem. Let's just mention that our heat kernel equals the fundamental solution in Q. Zhang's paper (since $G(x,t;y,s)$ satisfies the conjugate heat equation in the $(y,s)$ variables) and our result is an improvement, since there are no conditions on the Ricci curvature or the injectivity radius.  

Recently, Q. Zhang and X. Cao characterized in \cite{XCQZ10} the Type I singularity model of the Ricci flow by means of upper and lower bounds of the fundamental solution of the conjugate heat equation.  

In this paper, we obtain a bound on the heat kernel, depending on the best constants in a Sobolev imbedding theorem. We state the Sobolev imbedding theorems and the constants that the result will depend on in section (\ref{doi}), while the proof is given in section (\ref{trei}). We will conclude the paper with the special case, when the scalar curvature is positive at the initial time (and hence, since the Ricci flow preserves the positivity of the scalar curvature, throughout the manifold at all given times). 

We should add that the bound we get is not sharp, and that a long term goal is to get estimates similar to the Li-Yau ones. This, however, has proved to be much more difficult, due to the changing nature of the metric. Hope comes from the fact that similar gradient estimates on the solution of the heat equation have been found by the author, together with X. Cao and A. Pulemotov, using similar methods (see \cite{MBXCAP09}).

The main result of this paper can be stated as follows:  

\begin{theorem}\label{theorem}
Let $M^n$ be a closed Riemannian manifold, with $n\geq 3$ and let $\big(M,g(x,t)\big), t\in[0,T]$ be a solution to the Ricci flow~\eqref{Riccifloweqn}. Let $G(x,t;y,s)$ be the heat kernel, i.e.  fundamental solution for the heat equation~\eqref{heateqn}. Then there exists a positive number $C_n$, which depends only on the dimension $n$ of the manifold such that:
\begin{align*}
G(x,t;y,s)\leq &
\frac{C_n}{\left( \int\limits_{s}^{\frac{s+t}{2}}\left(\frac{m_0-c_n\tau}{m_0}\right)^{-2} \frac{e^{\frac{2}{n}H(\tau)}}{A(\tau)} \ d\tau\right)^{\frac{n}{4}} \left(\int\limits_{\frac{s+t}{2}}^{t}
\frac{e^{-\frac{2}{n}H(\tau)}}{A(\tau)} \ d\tau\right)^{\frac{n}{4}}} 
\end{align*}
for $0\leq s<t\leq T$; here $H(t)=\int\limits_{s}^{t}\left[\frac{B(\tau)}{A(\tau)}-\frac{3}{4}\cdot\frac{1}{m_0-c_n\tau}\right] d\tau $,  $1/m_0=\inf_{t=0}R$ - the infimum of the scalar curvature, taken at time $0$, and $A(t)$ and $B(t)$ are positive functions, which depend on the best constant in the Sobolev imbedding theorem.
\end{theorem}

One can notice that there are no curvature assumptions, just like in \cite{QZ06} and in \cite{SH08}, where the conjugate heat equation was analysed. 

The estimate may not seem natural, but in a special case, when the scalar curvature $R(x,0)>0$ (and thus $R(x,t)>0$ for any $t\in [0,T]$), one obtains a similar result to the one in the fixed metric case. Let's recall that J. Wang obtained in \cite{JW97} that the heat kernel on an $n$-dimensional compact Riemannian manifold $M$, with fixed metric, is bounded from above by $C(S)(t-s)^{-n/2}$, where $C(S)$ is the Neumann Sobolev constant of $M$, coming from a Sobolev imbedding theorem. Our corollary exibits a similar bound:   

\begin{corollary}
Under the same assumptions as in theorem (\ref{theorem}), together with the condition that the scalar curvature $R(x,t)$  be (strictly) positive at $t=0$, there exists a positive number $\tilde{C}_n$, which depends only on the dimension $n$ of the manifold and on the best constant in the Sobolev imbedding theorem in $\mathbb{R}^n$, such that:
\begin{align*}
G(x,t;y,s)\leq \tilde{C}_n\cdot\frac{1}{(t-s)^{\frac{n}{2}}} \hspace{1cm} \text{ for } 0\leq s<t\leq T
\end{align*}
\end{corollary}

The exact expression of $\tilde{C}_n$ will be shown in the proof. $\tilde{C}_n$ differs from $C(S)$ (as it appears in \cite{JW97}).

\vspace{10pt}



\section{Setup}\label{doi}

We consider an $n$-dimensional manifold without boundary $M$, which is compact, connected, oriented and smooth.

For $T>0$, let $\big(M,g(x,t)\big), t\in[0,T]$ to be a solution to the Ricci flow
\begin{align*}
\frac\partial{\partial t} g(x,t)=-2\Ric(x,t)
\end{align*}

The interval that we consider $[0,T]$ is a subset of the interval of short-time existence, hence we won't deal with blow-ups. 

We say that a smooth positive function $u:M\times[0,T]\to\mathbb R$ satisfies the heat equation if the following holds 
\begin{align*}
\left(\Delta-\frac\partial{\partial t}\right)u(x,t)=0
\end{align*} 

Here, $\Delta$ stands for the Laplacian given by $g(x,t)$. We will use $\nabla$ and $|\cdot|$ to denote the gradient and the norm with respect to $g(x,t)$, respectively. We emphasize that $\Delta$, $\nabla$, and $|\cdot|$ all depend on~$t\in[0,T]$. $XY$ denotes the scalar product of the vectors $X$ and $Y$ with respect to the metric $g(x,t)$.

We will denote the heat kernel, i.e. the fundamental solution of the heat operator $\left(\Delta-\frac\partial{\partial t}\right)$ by $G(x,t;y,s)$. Let's recall that the fundamental solution of an operator $L$ is a smooth function $G(x,t;y,s):M\times[0,T]\times M\times[0,T]\to\mathbb{R}$, with $s<t$, which satisfies two properties: 
\begin{enumerate}
\item[(i)] $L(G)=0$ in $(x,t)$ for $(x,t)\neq (y,s)$ 
\item[(ii)] $\lim_{t\to s}G(.,t;y,s)=\delta_y$ for every $y$, where $\delta_y$ is the Dirac delta function.
\end{enumerate} 

Guenther proved the existence and studied the properties of the fundamental solution for the operator $L(u)=\left(\Delta-\frac\partial{\partial t}-h(x,t)\right)u$ on a compact manifold whose metric evolves under the Ricci flow ($h(x,t)$ is a smooth function) \cite[Theorem 2.1]{CG02}. In particular, if $h(x,t)=0$ we get the existence of the heat kernel. 

During the computations, we will drop the arguments $(x,t;y,s)$ since it will be clear with respect to which variables we are considering the measure over which we are integrating.

Our proof relies on two Sobolev imbedding theorems, which are stated below. Since the manifold is compact, there will be no assumption on the injectivity radius or on the Ricci curvature (as in \cite{QZ06}). 

In \cite{TA76}, Aubin proved the following: 
\begin{theorem}\label{thm_Aubin}
Let $M^n$ be a compact Riemannian manifold. If $1\leq q\leq n$, then for any $\epsilon>0$, for any $q'\in[1,q]$ and for any $r>1$, there exists a positive constant $B(g,\epsilon, n)$ such that for any $\psi\in W^{1,q}(M)$ (the Sobolev space of weakly differentiable functions) :
\[||\psi||_p^r\leq  \left(K(n,q)^r+\epsilon\right)||\nabla\psi||_q^r+B||\psi||_{q'}^r \hspace{0.5cm}.\]

Here $K(n,q)$ is the best constant in the Sobolev imbedding (inequality) in $\mathbb{R}^n$ and $p=(nq)/(n-q)$.
\end{theorem}

Along the Ricci flow, Zhang proved the following uniform Sobolev inequality in \cite{QZ09}:

\begin{theorem}\label{thm_Zhang}
Let $M^n$ be a compact Riemannian manifold, with $n\geq 3$ and let $\big(M,g(t)\big)_{t\in[0,T]}$ be a solution to the Ricci flow~\eqref{Riccifloweqn}. Let $A$ and $B$ be positive numbers such that for $(M,g(0))$ the following Sobolev inequality holds: for any $v\in W^{1,2}(M,g(0))$,
\[\left(\int_M|v|^{\frac{2n}{n-2}}\ d\mu(g(0))\right)^{\frac{n-2}{n}}\leq A\int_M|\nabla v|^2\ d\mu(g(0))+B\int_Mv^2\ d\mu(g(0))\]

Then there exist positive functions $A(t)$, $B(t)$ depending only on the initial metric $g(0)$ in terms of $A$ and $B$, and $t$ such that, for all $v\in W^{1,2}(M,g(t))$, $t>0$, the following holds
 
\[\left(\int_M|v|^{\frac{2n}{n-2}}\ d\mu(g(t))\right)^{\frac{n-2}{n}}\leq A(t)\int_M\left(|\nabla v|^2+\frac{1}{4}Rv^2\right)\ d\mu(g(t))+B(t)\int_Mv^2\ d\mu(g(t))\]

Here $R$ is the scalar curvature with respect to $g(t)$. Moreover, if $R(x,0)>0$, then $A(t)$ and $B(t)$ are independent of $t$.
\end{theorem}

\section{Proof}\label{trei}


We assume without loss of generality that $s=0$. By the semigroup property of the heat kernel \cite[Theorem 2.6]{CG02} and the Cauchy-Bunyakovsky-Schwarz inequality we have that:
\begin{align*}
G(x,t;y,0)& =\int_M G\left(x,t;z,\frac{t}{2}\right)G\left(z,\frac{t}{2};y,0\right)d\mu\left(z,\frac{t}{2}\right) \\
          & \leq \left[\int_M G^2\left(x,t;z,\frac{t}{2}\right)d\mu\left(z,\frac{t}{2}\right) \right]^{1/2}\left[\int_M G^2\left(z,\frac{t}{2};y,0\right)d\mu\left(z,\frac{t}{2}
\right)\right]^{1/2}
\end{align*}

The key of the proof consists in determining upper bounds for the following two quantities:
\begin{align*}
\alpha(t)=\int_M G^2(x,t;y,s)d\mu(x,t)& \text{ (for $s$ fixed)}\\
\beta(s)=\int_M G^2(x,t;y,s)d\mu(y,s)& \text{ (for $t$ fixed)}
\end{align*}

Once we have these bounds, the conclusion follows immediately. The strategy to get these bounds consists in finding an ordinary differential inequality for each of the two quantities.

First let's recall that, by definition, $G$ satisfies the heat equation in the $(x,t)$ coordinates \[\Delta_x G(x,t;y,s)-\partial_t G(x,t;y,s)=0\] whereas in the $(y,s)$ it satisfies the conjugate heat equation \[\Delta_y G(x,t;y,s)+\partial_s G(x,t;y,s)- R(y,s) G(x,t;y,s)=0\] here $R(y,s)$ is the scalar curvature, measured with respect to the metric $g(s)$.

We will first deduce a bound on $\alpha(t)$, by finding an inequality involving $\alpha'(t)$ and $\alpha(t)$. Note that we will treat $G$ as being a function of $(x,t)$, the $(y,s)$ part is fixed. 

Since $\frac{d}{dt}(d\mu)=-Rd\mu$ (due to the Ricci flow), one has:

\begin{align}\label{pprim}
\alpha'(t) & =2\int_M G\cdot G_t\ d\mu(x,t)-\int_M G^2 R\ d\mu(x,t)= 2\int_M G\cdot(\Delta G)\ d\mu(x,t)-\int_M G^2 R\ d\mu(x,t)=\\
      & -2\int_M|\nabla G|^2\ d\mu-\int_M G^2R\ d\mu \leq -\int_M[|\nabla G|^2+R G^2]\ d\mu(x,t)
\end{align}

The difficult part will be to estimate $\int_M|\nabla G|^2d\mu$. The way we proceed is to use the Sobolev imbedding theorem, which gives a relation between $\int_M|\nabla G|^2d\mu$ and $\int_M G^2d\mu$, and the H\"older inequality to bound the term involving $G^{2n/(n-2)}$:

\begin{align}\label{Holder}
\int_M G^2\ d\mu(x,t)\leq \left[\int_M G^{\frac{2n}{n-2}}\ d\mu(x,t) \right]^{\frac{n-2}{n+2}} \left[\int_M G\ d\mu(x,t) \right]^{{4\over n+2}}
\end{align}

By theorem (\ref{thm_Aubin}), for $r=2$, $q=q'=2$ and $p=2n/(n-2)$ one gets that at time $t=0$, the following inequality holds for any $v\in W^{1,2}(M,g(0))$ and for some $B>0$:

\[\left(\int_M|v|^{\frac{2n}{n-2}}\ d\mu(g(0))\right)^{\frac{n-2}{n}}\leq (K(n,2)^2+\epsilon)\int_M|\nabla v|^2\ d\mu(g(0))+B\int_Mv^2\ d\mu(g(0))\]

where $K(n,2)$ - the best constant in the Sobolev imbedding in $\mathbb{R}$.

From this, by theorem \ref{thm_Zhang} it follows that at any time $t>0$ within the interval of short existance of the solution to the Ricci flow $[0,T)$ one has that for all $v\in W^{1,2}(M,g(t))$:

\[\left(\int_M|v|^{\frac{2n}{n-2}}\ d\mu(g(t))\right)^{\frac{n-2}{n}}\leq A(t)\int_M\left(|\nabla v|^2+\frac{1}{4}Rv^2\right)\ d\mu(g(t))+B(t)\int_Mv^2\ d\mu(g(t))\]

where $A(t)$ is a positive function depending on $g(0)$ and $(K(n,2)^2+\epsilon)$, while  $B(t)$ is also a positive function, depending on $B$.

Since $G(x,t;.,.)\in  W^{1,2}(M,g(t))$ (it is even smooth), then the above holds, so one can relate the RHS of (\ref{Holder}) to the Sobolev inequality:
\begin{align}\label{Sobolev1}
\int_M G^2\ d\mu(x,t) &  \leq \left[\int_M G^{\frac{2n}{n-2}}\ d\mu(x,t) \right]^{\frac{n-2}{n+2}} \left[\int_M G\ d\mu(x,t) \right]^{{4\over n+2}}\\
  & \leq \left[ A(t)\int_M\left(|\nabla G|^2+\frac{1}{4}RG^2\right)\ d\mu(x,t)+B(t)\int_MG^2\ d\mu(x,t)\right]^{\frac{n}{n+2}}\left[\int_M G\ d\mu(x,t)\right]^{\frac{4}{n+2}}
\end{align}

We need to estimate the term $J(t):=\int_M G(x,t;y,s)\ d\mu(x,t)$. By the definition of the fundamental solution, we have that: $\int_M G(x,t;y,s)d\mu(y,s)=1$, but that's not true if one integrates in $(x,t)$. We will obtain a differential inequality for $J(t)$ and the estimate will follow therefrom. 

\begin{align*}
J'(t)&=\int_M G_t(x,t;y,s)\ d\mu(x,t)+\int_M G(x,t;y,s)\frac{d}{dt}\ d\mu(x,t)=\int_M\Delta_xG(x,t;y,s)\ d\mu(x,t)\\
     &-\int_M G(x,t;y,s)R(x,t)\ d\mu(x,t)=-\int_M G(x,t;y,s)R(x,t)\ d\mu(x,t)
\end{align*}
the first term being $0$, as we are on a compact manifold, without boundary.

The scalar curvature satisfies the following differential inequality (see \cite{BCPLLN06}):
\begin{align*}
\frac{\partial R}{\partial t}-\Delta R-\frac{2}{n}R^2\leq 0
\end{align*}
Since the solutions of the ODE $\frac{d\rho}{dt}=\frac{2}{n}\rho^2$
are $\rho(t)=\frac{n}{n\rho(0)^{-1}-2t}$, by the maximum principle
we get a bound on the scalar curvature, for $s\leq \tau\leq t$:
\begin{align*}
R(z,\tau)\geq \frac{n}{n(\inf_{t=0}R)^{-1}-2\tau}=\frac{1}{(\inf_{t=0}R)^{-1}-\frac{2}{n}\tau}:=\frac{1}{m_0-c_n\tau}
\end{align*}
(here and later, if $\inf_{t=0} R\geq 0$, then the above is regarded
as zero).

\bigskip
Using this lower bound for $R$ (for $\tau\in(s,t]$), we get:
\begin{align*}
J'(\tau)\leq -\frac{1}{m_0-c_n\tau}J(\tau)
\end{align*}
After integrating the above from $s$ to $t$, while noting that by $J(s)$ one understands: 
\[J(s)=\lim\limits_{t\to s}\int_M G(x,t;y,s)\ d\mu(x,t)=\int_M \lim\limits_{t\to s} G(x,t;y,s)\ d\mu(x,t)=\int_M \delta_{y}(x)\ d\mu(x,s)=1\]
one obtains:
\begin{align*}
J(t)\leq \left(\frac{m_0-c_nt}{m_0-c_ns}\right)^{\frac{n}{2}}:=(\chi_{t,s})^{\frac{n}{2}}
\end{align*}

Hence $\int_M G(x,t;y,s)\ d\mu(x,t)\leq (\chi_{t,s}))^{\frac{n}{2}}$ and (\ref{Sobolev1}) becomes:

\begin{align*}
\int_M G^2d\mu(x,t)\leq \left[ A(t)\int_M\left(|\nabla G|^2+\frac{1}{4}RG^2\right)\ d\mu(x,t)+B(t)\int_MG^2\ d\mu(x,t)\right]^{\frac{n}{n+2}}\left(\chi_{t,s}\right)^{\frac{2n}{n+2}}
\end{align*}

From this it follows immediately that:

\begin{align*}
\int_M |\nabla G|^2\ d\mu(x,t)\geq \frac{1}{\chi^2_{t,s}A(t)}\left[\int_M G^2\ d\mu(x,t)\right]^{\frac{n+2}{n}} -\frac{B(t)}{A(t)}\int_M G^2\ d\mu(x,t)-\frac{1}{4}\int RG^2\ d\mu(x,t)
\end{align*}

Combining this with the inequality from (\ref{pprim}), one obtains the following differential inequality for $\alpha(t)$:
\begin{align*}
\alpha'(t)\leq -\frac{1}{\chi^2_{t,s}A(t)} \alpha(t)^{{n+2\over n}}+\frac{B(t)}{A(t)}\alpha(t)-\frac{3}{4}\int RG^2d\mu(x,t)
\end{align*}

Note that the above is true for any $\tau\in(s,t]$. For the following computation, we will consider $t$ fixed as well.  Recall that for $\tau\in (s,t]$, $R(\cdot, \tau)\geq\frac{1}{m_0-c_n\tau}$. Denoting with:
\begin{align*}
h(\tau):=\frac{B(\tau)}{A(\tau)}-\frac{3}{4}\cdot\frac{1}{m_0-c_n\tau}
\end{align*}
we get:
\begin{align*}
\alpha'(\tau)\leq -\frac{1}{\chi^{2}_{\tau,s}A(\tau)} \alpha(\tau)^{{n+2\over n}}+h(\tau)\alpha(\tau)
\end{align*}

Let $H(\tau)$ be an antiderivative of $h(\tau)$. By the integrating factor method, one finds:
\begin{align*}
(e^{-H(\tau)}\alpha(\tau))'\leq  -\frac{1}{\chi^{2}(\tau)A(\tau)}(e^{-H(\tau)}\alpha(\tau))^{\frac{n+2}{n}}e^{\frac{2}{n}H(\tau)}
\end{align*}
Since the above is true for any $\tau\in(s,t]$, by integrating from $s$ to $t$ and taking into account that 
\[\lim\limits_{\tau\searrow s}\alpha(\tau)=\int_M\lim\limits_{\tau\searrow s}G^2(x,\tau;y,s)\ d\mu(x,\tau)=\int_M\delta^2_y(x)\ d\mu(x,s)=0\]
one obtains the first necessary bound:
\begin{align*}
\alpha(t)\leq \frac{C_ne^{H(t)}}{\left(\int\limits_{s}^{t}\frac{e^{\frac{2}{n}H(\tau)}}{\chi^{2}(\tau)A(\tau)} d\tau\right)^{\frac{n}{2}}}
\end{align*}
where $C_n=\left(\frac{2}{n}\right)^{\frac{n}{2}}$.

\bigskip
The next step is to estimate $\beta(s)=\int_M G^2(x,t;y,s)\ d\mu(y,s)$, for which the computation is different, due to the assymetry of the
equation. As stated above, the second entries of $G$ satisfy the conjugated equation:
\begin{align*}
\Delta_y G(x,t;y,s)+\partial_s G(x,t;y,s)- R G(x,t;y,s)=0
\end{align*}

Proceeding just as in the $\alpha(t)$ case, we get the following:
\begin{align*}
\beta'(s)& =2\int_M GG_s\ d\mu(y,s)-\int_M RG^2\ d\mu(y,s) =2\int_M G(-\Delta G +RG)\ d\mu(y,s)-\int_M RG^2\ d\mu(y,s)\\
     & =-2\int_M G(\Delta G)\ d\mu(y,s)+\int_M RG^2\ d\mu(y,s)=2\int_M |\nabla G|^2\ d\mu(y,s)+\int_M RG^2\ d\mu(y,s) \\
     & \geq \int_M |\nabla G|^2\ d\mu(y,s)+\int_M RG^2\ d\mu(y,s)
\end{align*}

Hence
\begin{align*}
\beta'(s)\geq \int_M (|\nabla G|^2+RG^2)\ d\mu(y,s)
\end{align*}

But this time, by the property of the heat kernel:
\begin{align*}
\tilde{J}(s):=\int_M G(x,t;y,s)\ d\mu(y,s)=1
\end{align*}

so by applying H\"older (as for $\alpha(t)$) and relating it to the Sobolev inequality, we get:
\begin{align*}
\int_M G^2\ d\mu(y,s)&\leq \left[ A(s)\int_M\left(|\nabla G|^2+\frac{1}{4}RG^2\right)\ d\mu(y,s)+B(s)\int_MG^2\ d\mu(y,s)\right]^{\frac{n}{n+2}}\left[\int_M G\ d\mu(y,s)\right]^{\frac{4}{n+2}} \\
 & =\left[ A(s)\int_M\left(|\nabla G|^2+\frac{1}{4}RG^2\right)\ d\mu(y,s)+B(s)\int_MG^2\ d\mu(y,s)\right]^{\frac{n}{n+2}}
\end{align*}

Following the same steps as for $\alpha(t)$, one finds
\begin{align*}
\beta'(s)\geq \frac{1}{A(s)} \beta(s)^{{n+2\over n}}-h(s)\beta(s)
\end{align*}
($h(s)$ denotes, as before, $\frac{B(s)}{A(s)}-\frac{3}{4}\cdot\frac{1}{m_0-c_ns}$)

The above is true for any $\tau\in[s,t)$. We will apply again the integrating factor method, with $H(\tau)$ being the same antiderivative of $h(\tau)$ as above. 
For $\tau\in [s,t)$, the following holds:
\begin{align*}
(e^{H(\tau)}\beta(\tau))'\geq \frac{1}{A(\tau)}(e^{H(\tau)}\beta(\tau))^{{n+2\over n}}e^{-\frac{2}{n}H(\tau)}
\end{align*}

Integrating between $s$ and $t$, and taking into account that
\[ \lim\limits_{\tau\nearrow t}\beta(\tau)=\int_M\lim\limits_{\tau\nearrow t}G^2(x,t;y,\tau)\ d\mu(y,\tau)=\int_M\delta^2_y(x)\ d\mu(y,t)=0\] 
we get the second desired bound:
\begin{align*}
\beta(s)\leq \frac{C_ne^{-H(s)}}{\left(\int\limits_{s}^{t} \frac{e^{-\frac{2}{n}H(\tau)}}{A(\tau)}\ d\tau\right)^{n/2}}
\end{align*}

From the estimates of $\alpha$ and $\beta$ we get the following:

\begin{align*}
\alpha\left(\frac{t}{2}\right)=\int_M G^2\left(z,\frac{t}{2};y,0\right)\ d\mu\left(z,t/2\right)\leq \frac{C_ne^{H(t/2)}}{\left( \int\limits_{0}^{t/2}\left(\frac{m_0-c_n\tau}{m_0}\right)^{-2} \frac{e^{\frac{2}{n}H(\tau)}}{A(\tau)} \ d\tau\right)^{\frac{n}{2}}}
\end{align*}
\begin{align*}
\beta\left(\frac{t}{2}\right)=\int_M G^2\left(x,t;z,\frac{t}{2}\right)\ d\mu\left(z,\frac{t}{2}\right)\leq \frac{C_ne^{-H(t/2)}}{\left(\int\limits_{t/2}^{t}
\frac{e^{-\frac{2}{n}H(\tau)}}{A(\tau)} \ d\tau\right)^{n/2}}
\end{align*}

Here, we may choose $H(t/2)=\int\limits_{0}^{t/2}\left[\frac{B(\tau)}{A(\tau)}-\frac{3}{4}\cdot\frac{1}{m_0-c_n\tau}\right]\ d\tau$, since the relation is true for any antiderivative of $h(\tau)=\frac{B(\tau)}{A(\tau)}-\frac{3}{4}\cdot\frac{1}{m_0-c_n\tau}$.

The conclusion follows from multiplying the relations above.


\section{Special case: positive scalar curvature}

In the special case when $R(x,t)>0$, one gets that $J'(\tau)\leq 0$, which means that $J(\tau)$ is decreasing, so $J(\tau)\leq J(s)=1$, thus leading to the differential inequality for $\alpha(t)$ to be: 

\begin{align*}
\alpha'(t)\leq -\frac{1}{A(t)} \alpha(t)^{{n+2\over n}}+\frac{B(t)}{A(t)}\alpha(t)
\end{align*}

And from this the bound for $\alpha(t)$ becomes:

\begin{align*}
\alpha(t)\leq \frac{C_ne^{H(t)}}{\left(\int\limits_{s}^{t}\frac{e^{\frac{2}{n}H(\tau)}}{A(\tau)} d\tau\right)^{\frac{n}{2}}}
\end{align*}

where $H(\tau)$ is the antiderivative of $\frac{B(\tau)}{A(\tau)}$ such that $H(s)\neq 0$ and $H(t)\neq 0$.

Similarly, one obtaines for $\beta(s)$:

\begin{align*}
\beta'(s)\geq \frac{1}{A(s)} \beta(s)^{{n+2\over n}}-\frac{B(s)}{A(s)}\alpha(s)
\end{align*}

and from this: 

\begin{align*}
\beta(s)\leq \frac{C_ne^{-H(s)}}{\left(\int\limits_{s}^{t} \frac{e^{-\frac{2}{n}H(\tau)}}{A(\tau)}\ d\tau\right)^{n/2}}
\end{align*}

where $H(\tau)$ is the same antiderivative of $\frac{B(\tau)}{A(\tau)}$ as above. 

By (\ref{thm_Zhang}), in the case of $R(x,0)>0$, the two functions $A(t)$ and $B(t)$ are constants, let's call them $A$ and $B$. Recall that $A$ is in fact $K(n,2)+\epsilon$, where $K(n,2)$ is the best constant in the Sobolev imbedding and $\epsilon>0$.  

One has that $H(t)=\frac{B}{A}t$. Using this, we get:
\begin{align*}
G(x,t;y,s)\leq &
\frac{C_n}{\left( \int\limits_{s}^{\frac{s+t}{2}}\frac{e^{\frac{2}{n}H(\tau)}}{A(\tau)} \ d\tau\right)^{\frac{n}{4}} \left(\int\limits_{\frac{s+t}{2}}^{t}
\frac{e^{-\frac{2}{n}H(\tau)}}{A(\tau)} \ d\tau\right)^{\frac{n}{4}}}=\frac{C_n}{\left[\frac{n^2}{4B^2}\left(1-e^{-\frac{2B}{nA}\frac{t-s}{2}}\right)^2\right]^{\frac{n}{4}}}
\end{align*}

But, by Taylor expansion, the last expression is bounded by 
\[ \frac{C_n}{\left[\frac{n^2}{4B^2}\left(1-e^{-\frac{2B}{nA}\frac{t-s}{2}}\right)^2\right]^{\frac{n}{4}}}\leq \frac{\tilde{C}_n}{(t-s)^{\frac{n}{2}}}\]
where $\tilde{C}_n=C_n\cdot (2A)^{\frac{n}{2}}=\left(\frac{2}{n}\right)^{\frac{n}{2}}\cdot (2(K(n,2)+\epsilon))^{\frac{n}{2}}$.

Combining the two, we ge the desired corollary.


\bibliographystyle{alpha}
\bibliography{Fund_sol}

\end{document}